\documentclass[sn-mathphys,Numbered]{sn-jnl}

\usepackage{graphicx}%
\usepackage{multirow}%
\usepackage{amsmath,amssymb,amsfonts}%
\usepackage{amsthm}%
\usepackage{mathrsfs}%
\usepackage[title]{appendix}%
\usepackage{xcolor}%
\usepackage{textcomp}%
\usepackage{manyfoot}%
\usepackage{booktabs}%
\usepackage{algorithm}%
\usepackage{algorithmicx}%
\usepackage{algpseudocode}%
\usepackage{listings}%
\usepackage{mathtools}
\usepackage{subcaption}
\usepackage{todonotes}
\newcommand{\R}{\mathbb{R}}
\newcommand{\C}{\mathbb{C}}
\newcommand{\Mat}[1]{\begin{pmatrix*}[r]#1\end{pmatrix*}}
\renewcommand{\Vec}[1]{\begin{pmatrix*}[r]#1\end{pmatrix*}}
\renewcommand{\vec}[1]{\mathbf{#1}}
\newcommand{\mat}[1]{\mathbf{#1}}
\newcommand{\abs}[1]{\lvert#1\rvert}
\newcommand{\norm}[1]{\lVert#1\rVert}
\DeclareMathOperator{\diag}{diag}

\theoremstyle{thmstyleone}%
\theoremstyle{thmstyletwo}%

\theoremstyle{thmstylethree}%

\begin{document}

\title[]{On the non-global linear stability and spurious fixed points of MPRK schemes with negative RK parameters}


\author[1]{\fnm{Thomas} \sur{Izgin}}\email{izgin@mathematik.uni-kassel.de}

\author*[1]{\fnm{Stefan} \sur{Kopecz}}\email{kopecz@mathematik.uni-kassel.de}

\author[1]{\fnm{Andreas} \sur{Meister}}\email{meister@mathematik.uni-kassel.de}

\author[1]{\fnm{Amandine}\sur{Schilling}}\email{schilling@uni-kassel.de}

\affil*[1]{\orgdiv{Institute of Mathematics}, \orgname{University of Kassel}, \orgaddress{\street{Untere K\"onigsstr. 86}, \postcode{34117} \city{Kassel}, \country{Germany}}}


\abstract{Recently, a stability theory has been developed to study the linear stability of modified Patankar--Runge--Kutta (MPRK) schemes. This stability theory provides sufficient conditions for a fixed point of an MPRK scheme to be stable as well as for the convergence of an MPRK scheme towards the steady state of the corresponding initial value problem, whereas the main assumption is that the initial value is sufficiently close to the steady state. Initially, numerical experiments in several publications indicated that these linear stability properties are not only local, but even global, as is the case for general linear methods. Recently, however, it was discovered that the linear stability of the MPDeC(8) scheme is indeed only local in nature. 
Our conjecture is that this is a result of negative Runge--Kutta (RK) parameters of MPDeC(8) and that linear stability is indeed global, if the RK parameters are nonnegative.
To support this conjecture, we examine the family of MPRK22($\alpha$) methods with negative RK parameters and show that even among these methods there are methods for which the stability properties are only local.  However, this local linear stability is not observed for MPRK22($\alpha$) schemes with nonnegative Runge-Kutta parameters.}

\keywords{Conservative scheme, unconditional positivity, modified Patankar--Runge--Kutta methods, linear stability}



\maketitle

\section{Introduction}
Modified Patankar--Runge--Kutta (MPRK) schemes are numerical time integration schemes which are conservative and unconditionally positivity-preserving, when applied to a positive and conservative production-destruction system (PDS)
\begin{equation}\label{eq:PDS}
 y_i'=\sum_{\substack{j=1\\j\ne i}}^N\bigl( p_{ij}(t,\mathbf y) - d_{ij}(t,\mathbf y)\bigr),\quad i=1,\dots, N,\end{equation}
where $\mathbf y=(y_1,\dots,y_N)^T$ and $p_{ij}(t,\vec y),d_{ij}(t,\vec y)\geq 0$ for all $\vec y>\vec 0$, $t\geq 0$.
A PDS is called positive, when positive initial values imply positive solutions for all times, and is called conservative, when $\mathbf 1^T\mathbf y'=0$ or equivalently $\mathbf 1^T\mathbf y$ is constant for all times. 
Here, $\mathbf 1=(1,\dots,1)^T\in\R^N$ denotes the vector with all elements equal to one. 
Since MPRK schemes do not belong to the class of general linear methods, they can be unconditionally positivity-preserving and of higher order at the same time. For general linear methods unconditional positivity is restricted to first order schemes, see \cite{bolleycrouzeix78,Sandu2002,HundsdorferVerwer2003}. So far, second and third order MPRK schemes for the integration of autonomous PDS have been constructed  in \cite{BDM2003,KopeczMeister2018,KopeczMeister2018b}, based on the classical Runge--Kutta (RK) formulation.  The same was done in \cite{HZS2019,HuangShu2019}, based on the SSP formulation of RK schemes. In \cite{TorloOeffner2020} arbitrary high order MPRK schemes on the basis of deferred correction schemes have been introduced. MPRK schemes for time dependent PDS have been investigated in \cite{AGKM2021}.
For other approaches to obtain unconditional positivity we refer to \cite{Sandu2001,Sandu2002,BBKS2007,BRBM2008,MCD2020, AKM2020, BAM2022}. We also want to mention the novel framework \cite{NSARK} on order conditions for Runge--Kutta-like schemes to which MPRK methods belong. Besides a general theory for deriving order conditions, for the first time sufficient and necessary conditions for fourth order MPRK schemes were presented and reduced therein. Additionally, the authors proved that the order of an MPRK method does not depend on the signs of the RK parameters.

In addition to the order of a numerical scheme, its stability is of course crucial for its usefulness in practical applications. 
In the following, we will be concerned with linear stability of MPRK schemes, i.\,e. the stability behavior of MPRK schemes when applied to a positive and conservative linear system.
Based on $\mat A=(a_{ij})\in\R^{N\times N}$, a positive and conservative linear system has the form
\begin{equation}\label{eq:linsys}
\mathbf y'=\mathbf A\mathbf y,\quad \mathbf A - \diag(\mathbf A)\geq 0,\quad \mathbf  1^T\mathbf A=\mathbf 0,
\end{equation}
where $\mathbf A-\diag(\mathbf A)\geq 0$ is necessary and sufficient for positivity, see \cite{luenberger1979introduction}, and $\mathbf  1^T\mathbf A=\mathbf 0$ must hold to ensure conservativity. 
Here, the notation $\diag(\mathbf A)\in\R^{N\times N}$ is used to denote the diagonal matrix with diagonal elements equal to those of $\mathbf A$.
Furthermore, the conservation property $\mathbf  1^T\mathbf A=\mathbf 0$ together with positivity implies $\diag(\mathbf A)\leq 0$. 
If we define $\mathbf B=\mathbf A-\diag(\mathbf A)$, system \eqref{eq:linsys} can be rewritten in production-destruction form as
\begin{equation}\label{eq:linpds}
\mathbf y'=\mathbf B\mathbf y - (-\diag(\mathbf A)\mathbf y),\quad \mathbf B\geq 0, \quad\mathbf 1^T\mathbf A=\mathbf 0.
\end{equation}
Since $\mathbf B=(b_{ij})\geq \mathbf 0$ and $\mathbf y\geq \mathbf 0$ by assumption the system's production terms are $p_{ij}(\vec y)=b_{ij}y_j=a_{ij}y_j\geq 0$ for $i\ne j$ and the corresponding destruction terms $d_{ij}$ are contained within $-\diag(\mathbf A)\mathbf y\geq \mathbf 0$. For $N=2$ all positive and conservative linear systems \eqref{eq:linpds} can be represented as
\begin{equation}\label{eq:lin2x2}
\vec y'= \Mat{-a & b\\ a & -b}\vec y,\quad a,b\geq 0.
\end{equation}

The central notion for linear stability of general linear methods is A-stability. A general linear method is said to be A-stable, whenever its numerical solution of $y'=\lambda y$ for an arbitrary time step size $\Delta t>0$ tends to zero for all $\lambda \in\C^-=\{z\in\C\mid \operatorname{Re}(z)<0\}$. 
The choice $\lambda\in\C^-$ comes from the fact, that A-stability ensures that the numerical solution of a linear system $\vec y'=\mat A\vec y$ tends to $\vec 0$, whenever all eigenvalues of $\mat A$ belong to $\C^-$. 
For an A-stable linear method applied to $\mathbf y'=\mathbf A\mathbf y$ with a spectrum $\sigma(\mathbf A)\subseteq \C^-$, the unique steady state solution $\vec y^*=\vec 0$ is an asymptotically stable fixed point for arbitrary time step sizes. 
However, the crucial difference between A-stability and asymptotic stability is that the former is a global property, i.e., independent of the initial value, while the latter is a local property, since asymptotic stability requires that the initial value be sufficiently close to the fixed point.

Due to its importance for general linear methods, we would also like to investigate the A-stability of MPRK schemes. 
Unfortunately, several obstacles stand in the way.
First, MPRK schemes cannot be applied to the scalar linear test equation, particularly not with $\lambda\in\C^-$, as it is unclear how the complex term $\lambda y$ can be split into production and destruction terms. 
To get around this, we can apply MPRK schemes directly to positive and conservative linear systems. But, the conservation property is in contradiction with the asymptotic stability of $\vec y^*=\vec 0$, i.\,e.\ there is no conservative linear system whose solution tend to $\vec 0$ for $t\to\infty$ for initial values $\vec y^0>\vec 0$, since $\vec 1^T\vec y$ must be constant for all times.
Moreover, a conservative linear system can possess several independent linear invariants, i.\,e.\ there exist $K$ linear independent vectors $\vec n_i$ with $\vec n_i^T\mat A=\vec 0$ and hence, $\vec n_i^T\vec y$ remains constant for all times. 
To avoid the issue with asymptotic stability, we can weaken the requirement and demand only stability instead, while at the same time, we can require that the numerical approximations tend to the unique steady state solution $\vec y^*$ of \[\vec y'=\mat A\vec y,\quad \mathbf y(0)=\mathbf y^0,\]
which always satisfies $\vec n_i^T\vec y^*=\vec n_i^T\vec y^0$. 
In summary, we are looking for MPRK schemes for which a stable steady state $\vec y^*$ of a positive and conservative linear system $\vec y'=\mat A\vec y$,  becomes a stable fixed point of the MPRK scheme for all time step sizes $\Delta t$, and in addition, the iterates $\vec y^n$ tend to $\vec y^*$ for all initial values $\vec y^0$ that satisfy $\vec n_i^T\vec y^0=\vec n_i^T\vec y^*$.
 
The fact that an MPRK scheme $\mathbf y^{n+1}=\mathbf g(\mathbf y^n)$ is not a general linear method makes it complicated to find MPRK schemes with the desired properties,
since the application of an MPRK scheme to a linear system results in a nonlinear iteration of the form
\[\mathbf y^{n+1}=\mathbf R(\Delta t\mathbf A,\mathbf y^n)\mathbf y^n.\]
In addition, the conservation property implies the existence of infinitely many non-hyperbolic fixed points $\mathbf y^*\ne\vec 0$ of the map $\mathbf g$, i.\,e.\ the Jacobian $\mathbf D\mathbf g(\mathbf y^*)$ has eigenvalues with absolute value equal to one. Hence, a linear stability theory for MPRK schemes must be a stability theory for non-hyperbolic fixed points of nonlinear iterations.  
One such approach to study stability is based on the center manifold theory of dynamical systems and was introduced in \cite{IKM2022,IKM2022b}. 
Assuming that the initial value $\vec y^0$ is sufficiently close to the fixed point $\vec y^*$, the theory provides sufficient conditions for the stability of $\vec y^*$, as well as for the convergence of the iterates to the steady state of the corresponding initial value problem. Thus, this stability theory would almost suffice to find the desired schemes if the theory did not make local statements only, with respect to the initial value.
However, to the authors' knowledge this is the only approach to study stability of MPRK schemes so far.

The stability theory of \cite{IKM2022,IKM2022b} was used to investigate the linear stability of MPRK22 schemes and it was proven therein that MPRK22($\alpha$) schemes with $\alpha\geq 0.5$ are linearly stable and that their iterates tend to the correct fixed point $\vec y^*$, whenever the initial value $\vec y^0$ is sufficiently close to $\vec y^*$.
The stability theory was used in \cite{HIKMS2023} to find SSP-MPRK schemes with the same properties.
A stability analysis of the third order MPRK schemes from \cite{KopeczMeister2018b} and the MPDeC schemes of \cite{TorloOeffner2020} was carried out in \cite{OeffnerIzgin2023}.
We also want to note, that the stability theory of \cite{IKM2022,IKM2022b} is not only applicable to MPRK schemes and was also used in \cite{IKMM2023} to analyze the linear stability of BBKS \cite{BBKS2007,BRBM2008,AGKM2021} and GeCo schemes \cite{MCD2020}. Furthermore, it was used in \cite{IKM2022c} to investigate the stability of an MPRK scheme in the context of a nonlinear PDS.

Even though the stability results of \cite{IKM2022,IKM2022b} are only valid in a sufficiently small neighborhood of the fixed point, the numerical results in \cite{IKM2022,IKM2022b,HIKMS2023} suggested that the local stability might actually be a global one, just like it is the case for A-stable general linear methods.

However, in \cite{OeffnerIzgin2023,TOR2022} it was discovered that the MPDeC(8) method with equidistant nodes is indeed only locally stable.
To check the conjecture that this is due to the negative RK parameters of the MPDeC(8) scheme, the family of MPRK22($\alpha$) schemes with negative RK parameters, i.\,e. $\alpha<0.5$, was investigated in \cite{Schilling2023}. 
The numerical experiments presented therein show that MPRK22($-0.5$) is another MPRK scheme, which is only locally stable.

The aim of this paper is to summarize and extend the results of \cite{Schilling2023} and to justify the conjecture that this local linear stability behavior, which is unknown from general linear methods, only occurs for MPRK schemes if the Butcher tableau of the underlying RK scheme contains negative values. 

\section{Linear stability of MPRK22($\alpha$) schemes}
This section summarizes some of the results of \cite{Schilling2023}.

The idea of MPRK schemes is to modify explicit RK schemes by 
introducing additional weighting factors that ensure unconditional positivity and conservation. Destruction terms are multiplied by weights with respect to the equation they appear in. Production terms are multiplied by the same weights as their corresponding destruction counter parts. 
If all parameters of the underlying RK schemes are nonnegative, there is no difference between production or destruction terms on the continuous and discrete level. But if production or destruction terms are multiplied by a negative RK parameter they switch their roles from the continuous to the discrete level, which has to be dealt with appropriately within the implementation of the scheme. 

The MPRK22$(\alpha)$ schemes were introduced in \cite{KopeczMeister2018} and are based on general second order explicit RK schemes, i.\,e. $a_{21}=\alpha\ne 0$, $b_2=\frac{1}{2\alpha}$, $b_1=1-b_2$.  In \cite{KopeczMeister2018} only nonnegative RK parameters were considered, which is the case, if $\alpha\geq\frac 12$. For $0<\alpha<\frac12$ the parameter $b_1$ becomes negative and for $\alpha<0$ the parameters $a_{21}$ and $b_1$ are negative. Consequently, we need to distinguish these three different cases, for varying values of $\alpha$.

In a form which is suitable for positive as well as negative RK parameters, the MPRK22($\alpha$) schemes for the solution of \eqref{eq:PDS} can be reformulated as
\begin{subequations}\label{eq:MPRK}
\begin{align}\label{eq:MPRKstages}
y_i^{(1)} &=y_i^n,\\
y_i^{(2)} &= y_i^n + a_{21}\Delta t \sum_{j=1}^N \biggl(p_{ij}(\mathbf y^{(1)})\frac{y_{\gamma(j,i,a_{21})}^{(2)}}{ y^{(1)}_{\gamma(j,i,a_{21})}}-d_{ij}(\mathbf y^{(1)})\frac{y_{\gamma(i,j,a_{21})}^{(2)}}{ y^{(1)}_{\gamma(i,j,a_{21})}}\biggr),\\\label{eq:MPRKapprox}
 y_i^{n+1}&=y_i^n+\Delta t \sum_{k=1}^2 b_k 
 \sum_{j=1}^N\biggl( p_{ij}(\mathbf y^{(k)})\frac{y_{\gamma(j,i,b_k)}^{n+1}}{\sigma_{\gamma(j,i,b_k)}}-d_{ij}(\mathbf y^{(k)})\frac{y_{\gamma(i,j,b_k)}^{n+1}}{\sigma_{\gamma(i,j,b_k)}}\biggr),
\end{align}
\end{subequations}
for $i=1,\dots,N$, with
\[\sigma_i = \sigma_i(\mathbf y^n,\mathbf y^{(2)})= (y_i^n)^{1-1/a_{21}}(y_i^{(2)})^{1/a_{21}},\quad i=1,\dots,N\]
and the index function
\begin{equation}\label{eq:Indexfunktion}
\gamma(i,j,\theta)=i \ \  \text{for } \ \theta \geq 0 \ \  \text{und } \ \gamma(i,j,\theta)=j \ \  \text{for } \ \theta < 0.
\end{equation}
The purpose of the index function is to decide, whether a term is a production or destruction term and to choose the weighting factors accordingly. We also note that the index function \eqref{eq:Indexfunktion} was introduced in \cite{TorloOeffner2020} for the definition of MPDeC schemes.

It was proven in \cite{KopeczMeister2018} that MPRK22$(\alpha)$ schemes with $\alpha\geq\frac12$ are unconditionally positive and conservative second order schemes. The same is true for $0<\alpha<\frac12$ and $\alpha<0$, which was proven in \cite{Schilling2023} by a straight-forward modification of the proof given in \cite{KopeczMeister2018}. For a general framework concerning the order conditions of MPRK schemes we refer to \cite{NSARK}.


To examine the linear stability of MPRK22$(\alpha)$ schemes in terms of the stability theory \cite{IKM2022,IKM2022b}, we consider their application to positive and conservative linear PDS \eqref{eq:linpds}.
In terms of a shorter notation we define
\[\diag(\mathbf v/\mathbf u)=\diag(\mathbf v)\diag(\mathbf u)^{-1},\]
where $\diag(\mathbf v)=\diag(v_1,\dots,v_n)\in\R^{N\times N}$ denotes the diagonal matrix with the elements of the vector $\mathbf v=(v_1,\dots,v_n)^T\in\R^N$ on the diagonal. As mentioned above, we need to distinguish between the following three cases.
\subsection{Case $\alpha\geq 0.5$}
This case was already considered in \cite{IKM2022,IKM2022b} and we present the results for the sake of completeness.

As $\alpha=a_{21}>0$ implies that all RK parameters are nonnegative, the application of an MPRK22($\alpha$) scheme \eqref{eq:MPRK} to \eqref{eq:linpds} yields
\begin{align*}
\mathbf y^{(2)}&=\mathbf y^n + a_{21}\Delta t\bigl(\mathbf B\diag(\mathbf y^{(2)}/\mathbf y^n)\mathbf y^n+\diag(\mathbf y^{(2)}/\mathbf y^n)\diag(\mathbf A)\mathbf y^n\bigr),\\
\mathbf y^{n+1}&\!\begin{multlined}[t][11cm]
     = \mathbf y^n+\Delta t\bigl(\mathbf B\diag(\mathbf y^{n+1}/\boldsymbol \sigma)(b_1 \mathbf y^n+b_2\mathbf y^{(2)})\\+\diag(\mathbf y^{n+1}/\boldsymbol \sigma)\diag(\mathbf A)(b_1 \mathbf y^n+b_2\mathbf y^{(2)})\bigr),
     \end{multlined}
\end{align*}
In comparison to the underlying RK scheme a diagonal matrix with Patankar-weights was introduced on the left of $\diag(\mathbf A)$ in the destruction parts and on the right of $\mathbf B$ in the production parts. This is done in the other cases as well.

The dependence of $\mathbf y^{n+1}$ on $\mathbf y^n$ can be expressed by an implicit function $\mathbf g$, i.\,e. $\mathbf y^{n+1}=\mathbf g(\mathbf y^n)$. 
Each steady state $\mathbf y^*$ of \eqref{eq:linpds} is also a fixed point of $\mathbf g$ and
\begin{equation*}
\mathbf D\mathbf g(\mathbf y^*)
= \bigl(\mathbf I - \Delta t \mathbf A \big )^{-1}
\bigl(\tfrac{1}{2\alpha}  \Delta t \mathbf A \bigl(\mathbf I
-\bigl(\mathbf I-\alpha \Delta t \mathbf A \bigr)^{-1}\bigr) + \mathbf I \bigr)
\end{equation*}
has the eigenvalues $R(\Delta t \lambda)$, where $\lambda$ is an eigenvalue of $\mathbf A$
and the stability function is
\begin{equation*}
R(z)
=
\frac{\tfrac{1}{2\alpha} z
\left(1-\tfrac{1}{1-\alpha z}\right) +1}
{1- z}
=
\frac{-z^2-2\alpha z+2}
{2(1-\alpha z)(1-z)}.
\end{equation*}
As a result, $\abs{R(z)}<1$ for all $z\in\C^-$, which implies stability of fixed points and convergence towards the steady state of the underlying initial value problem.
Hence, the MPRK22($\alpha$) schemes with $\alpha\geq 0.5$ fulfill all desired properties, apart from the fact that stability could only be proven in a sufficiently small neighborhood of the fixed point.
Furthermore, we want to emphasize that in this case, apart from $\Delta t$, stability depends only on $\lambda$ like in the continuous case.

\subsection{Case $0<\alpha<0.5$}\label{sec:0<alpha<0.5}
A similar analysis as for $\alpha \geq 0.5$ can be conducted for the situation that $0<\alpha<0.5$. This was done in \cite{Schilling2023} and the results will be summarized here. In this case we have $b_1<0$ and $a_{21},b_2>0$. Consequently, the scheme \eqref{eq:MPRK} applied to \eqref{eq:linpds} reads
\begin{align*}
\mathbf y^{(2)}&=\mathbf y^n + a_{21}\Delta t\bigl(\mathbf B\diag(\mathbf y^{(2)}/\mathbf y^n)\mathbf y^n+\diag(\mathbf y^{(2)}/\mathbf y^n)\diag(\mathbf A)\mathbf y^n\bigr),\\
\mathbf y^{n+1}&\!\begin{multlined}[t][11cm]
     = \mathbf y^n+\Delta t\bigl(b_1\diag(\mathbf y^{n+1}/\boldsymbol \sigma)\mathbf B\mathbf y^n+b_2\mathbf B\diag(\mathbf y^{n+1}/\boldsymbol \sigma)\mathbf y^{(2)}+\\b_1\diag(\mathbf A)\diag(\mathbf y^{n+1}/\boldsymbol \sigma)\mathbf y^n+b_2\diag(\mathbf y^{n+1}/\boldsymbol \sigma)\diag(\mathbf A)\mathbf y^{(2)})\bigr).
     \end{multlined}
\end{align*}
The stability behavior is more complicated than for $\alpha\geq 0.5$ as the Jacobian $\mathbf D\mathbf g(\mathbf y^*)$ in general also depends on $\mathbf A^T$ and the steady state $\mathbf y^*$ itself, see \cite{Schilling2023} for details.
However, if we restrict ourselves to the specific system \eqref{eq:lin2x2}, then 
\begin{multline*}
\mat D\vec g(\vec y^*)
= -\bigl(-\mat I + \Delta t \bigl( -1 + \tfrac{1}{\alpha} \bigr) \mat A \big )^{-1}
\bigl(\mat I + \Delta t \bigl(2 - \tfrac{5}{2\alpha} + \tfrac{1}{\alpha^2} \bigr) \mat A \\
-\Delta t \bigl(\tfrac{3}{2\alpha} - \tfrac{1}{\alpha^2} \bigr) \mat A
\bigl(-\mat I+\alpha \Delta t \mat A \bigr)^{-1} \bigr)
\end{multline*}
and
\begin{align*}
R(z)=
-\frac{1+\bigl(2- \tfrac{5}{2\alpha} + \tfrac{1}{\alpha^2} \bigr)z
-\bigl( \tfrac{3}{2\alpha} - \tfrac{1}{\alpha^2} \bigr)
\frac{z}{(-1+\alpha z)} }
{-1+\bigl(-1+\frac{1}{\alpha}\bigr)z}.
\end{align*}
A technical computation shows that for every $\alpha$ there exists a $z^*<0$ for which $\abs{R(z^*)}>1$. Hence, the MPRK22($\alpha$) schemes are only conditionally stable for $0<\alpha< 0.5$.
To be precise, $z^*$ is given by
\[z^* = \frac{-2 \alpha^2 +3 \alpha - 2 - \sqrt{4 \alpha^4+ 12 \alpha^3 - 11 \alpha^2 - 4 \alpha +4}}
{6 \alpha^2 - 7 \alpha +2}.
\]
\subsection{Case $\alpha<0$}\label{sec:alpha_neg}
In this case we have $a_{21},b_2<0$ and $b_1>0$. Hence, terms multiplied by $a_{21}$ or $b_2$ change their role from the continuous to the discrete level. As a result, MPRK22($\alpha$) schemes \eqref{eq:MPRK} applied to \eqref{eq:linpds} become
\begin{subequations}\label{eq:MPRK22<0}
\begin{align}
\mathbf y^{(2)}&=\mathbf y^n + a_{21}\Delta t\bigl(\diag(\mathbf y^{(2)}/\mathbf y^n)\mathbf B\mathbf y^n+\diag(\mathbf A)\diag(\mathbf y^{(2)}/\mathbf y^n)\mathbf y^n\bigr),\\
\mathbf y^{n+1}&\!\begin{multlined}[t][11cm]
     = \mathbf y^n+\Delta t\bigl(b_1\mathbf B\diag(\mathbf y^{n+1}/\boldsymbol \sigma)\mathbf y^n+b_2\diag(\mathbf y^{n+1}/\boldsymbol \sigma)\mathbf B\mathbf y^{(2)}\\+b_1\diag(\mathbf y^{n+1}/\boldsymbol \sigma)\diag(\mathbf A)\mathbf y^n+b_2\diag(\mathbf A)\diag(\mathbf y^{n+1}/\boldsymbol \sigma)\mathbf y^{(2)}\bigr).
     \end{multlined}
\end{align}
\end{subequations}
Again we summarize the results of \cite{Schilling2023} here. But as this is the important case for the purpose of this paper, we also give proofs in Appendix~\ref{apndx}, at least for specific cases.

As for $0<\alpha<0.5$ the Jacobian $\mathbf D\mathbf g(\mathbf y^*)$ in general depends on $\mathbf A^T$ and the steady state $\mathbf y^*$ itself.
But for the purpose of this paper it is sufficient to restrict ourselves to system \eqref{eq:lin2x2}, in which case it can be seen that the dependence on $\mat A^T$ and $\vec y^*$ within \eqref{eq:Jacobians} disappears due to the property \eqref{eq:condition} and  we have 
\begin{multline*}
\mat D\vec g(\vec y^*)
= -\bigl(-\mat I + \Delta t \bigl( 1 - \tfrac{1}{\alpha} \bigr) \mat A \big )^{-1}
\bigl(\mat I + \Delta t \bigl( \tfrac{3}{2\alpha} - \tfrac{1}{\alpha^2} \bigr) \mat A \\
-\Delta t \bigl( -\tfrac{1}{2\alpha} + \tfrac{1}{\alpha^2} \bigr) \mat A
\bigl(-\mat I-\alpha \Delta t \mat A \bigr)^{-1} \bigl( \mat I+2\alpha \Delta t \mat A \bigr) \bigr)
\end{multline*}
as well as
\begin{equation*} 
 R(z)=
 -\frac{1+\bigl( \tfrac{3}{2\alpha} - \frac{1}{\alpha^2} \bigr)z
 -\bigl( -\tfrac{1}{2\alpha} + \frac{1}{\alpha^2} \bigr)z
 \frac{1}{(-1-\alpha z)} \bigl(1+2\alpha z\bigr)}
 {-1+\bigl(1-\frac{1}{\alpha}\bigr)z}.
\end{equation*}
The nonzero eigenvalue of $\mathbf A$ in \eqref{eq:lin2x2} is $\lambda=-(a+b)<0$. Hence, we only need to consider $z\in\R$ with $z<0$. For every $\alpha<0$ the function $R$ is monotonically increasing on $(-\infty,0)$ with 
\[\lim_{z\to-\infty}R(z)=-\frac{\alpha+2}{2\alpha(\alpha-1)}\]
and consequently 
\begin{equation*}
-\frac{\alpha+2}{2\alpha(\alpha-1)}<R(z)<1 \quad \text{for } z<0 \text{ and } \alpha<0.
\end{equation*}
If we restrict $\alpha$ additionally to $\alpha\leq -0.5$, then $\abs{R(z)}<1$ for $z<0$ and it follows that MPRK22($\alpha$) schemes with $\alpha\leq-0.5$ are unconditionally stable, i.\,e. the stability is independent of $\Delta t$. 
Moreover, convergence to the steady state of the corresponding initial value problem is guaranteed as well, whenever the initial values is sufficiently close to the fixed point.
The situation is different for $-0.5<\alpha<0$, for which the schemes are only conditionally stable. For $-0.5<\alpha<0$ stability requires $z^*<z<0$ with
\[
z^*= \frac{2 \alpha^2 - \alpha +2+  \sqrt{4 \alpha^4+ 4 \alpha^3 - 3 \alpha^2 - 12 \alpha +4}}
{2 \alpha^2 - 3 \alpha -2}.
\]
As mentioned before, see \cite{Schilling2023} or Appendix~\ref{apndx} for details.

Now we know that MPRK22($\alpha$) schemes with $\alpha\leq -0.5$ are unconditionally stable and converge to the steady state of the corresponding initial value problem, at least when applied to the linear system \eqref{eq:lin2x2} and the initial value is not to far away from the steady state. 
As mentioned above, the numerical results in \cite{IKM2022,IKM2022b,HIKMS2023} suggested that this local behavior might actually be the global behavior.
But it was observed in \cite{Schilling2023} and \cite{OeffnerIzgin2023,TOR2022} that MPRK22(-0.5) and MPDeC(8) with equidistant nodes are indeed only locally stable. Such a linear stability behavior is unknown from general linear methods.

The conjecture that this can only occur in the presence of negative RK parameters is supported in Section~\ref{sec:numres}, where we also show that there are further values of $\alpha$ with $\alpha<-0.5$ for which the same stability behavior can be observed. Moreover, the experiments with $\alpha\geq 0.5$, corresponding to non-negative RK parameters, do not show a case for which the stability properties are local.
  
\section{Numerical results}\label{sec:numres}
In this section we demonstrate numerically that for some MPRK22($\alpha$) schemes with $\alpha<0$ the linear stability and convergence to the correct steady state is indeed only given locally, i.\,e.\ the initial value must be sufficiently close to the fixed point.
To show this, it is already sufficient to consider the linear system
\begin{subequations}\label{eq:testprob}
\begin{equation}\vec y' = \Mat{-a & a\\a & -a}\vec y,\quad a>0\end{equation}
with initial value
\begin{equation}\vec y^0=\frac12\Vec{1\\1}+\delta\Vec{1\\-1},\quad 0\leq\delta<\frac12,\end{equation}
\end{subequations}
where the restriction on $\delta$ is necessary to ensure positive initial data. 
The solution of this initial value problem is
\[\vec y(t)=\frac{1}{2}\Vec{1\\1}+\delta e^{-2 a t}\Vec{1\\-1}\] 
and its steady state is
\[\vec y^*=\lim_{t\to\infty} \vec y(t)=\frac12\Vec{1\\1}\]
independent of $\delta$.

The following computations were performed with MATLAB 2023b. The MATLAB code is available from \texttt{https://github.com/SKopecz/locstabMPRK22.git}.

\begin{figure}
\centering
\begin{subfigure}{0.5\textwidth}
    \includegraphics[width=\textwidth]{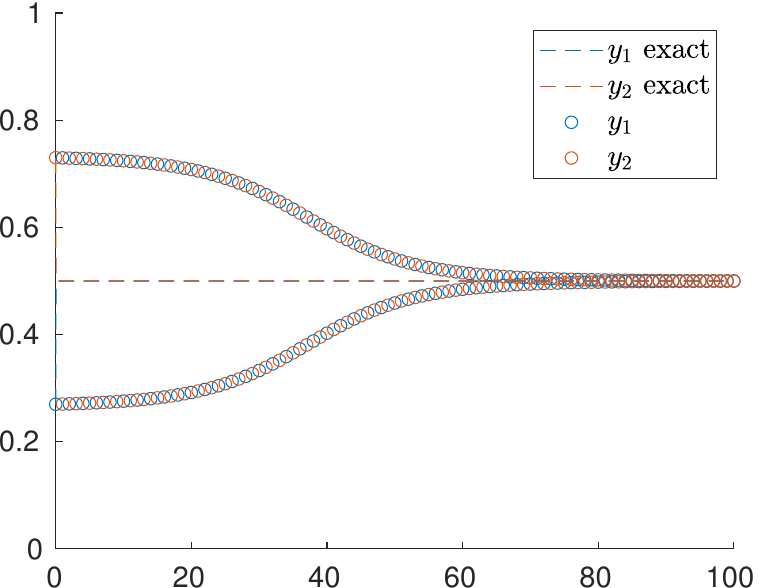}
    \caption{$\delta=0.23$}
    \label{fig:fig1a}
\end{subfigure}%
\begin{subfigure}{0.5\textwidth}
    \includegraphics[width=\textwidth]{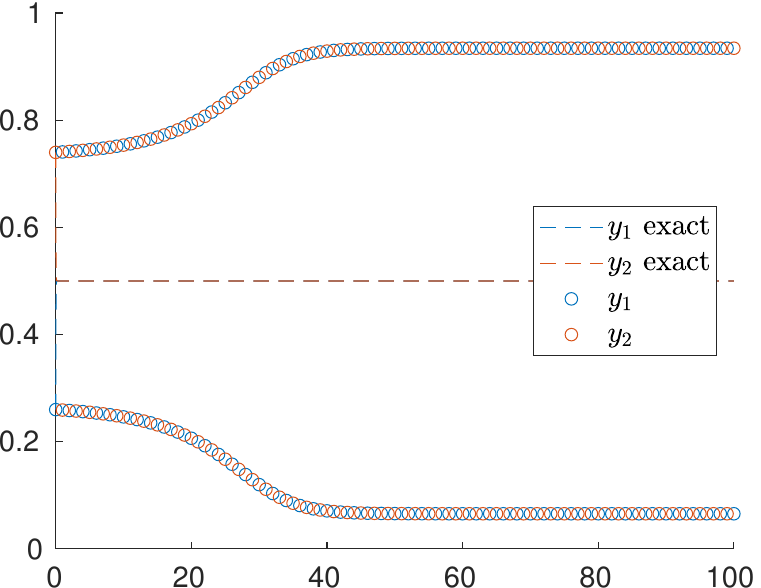}
    \caption{$\delta=0.24$}
    \label{fig:fig1b}
\end{subfigure}
\caption{Numerical solutions of the initial value problem \eqref{eq:testprob} and $a=20$ computed with MPRK22($-0.5$) and $\Delta t = 1$ for $\delta=0.23$ and $\delta=0.24$. }
\label{fig:fig1}
\end{figure}
First, we consider system \eqref{eq:testprob} with $a=20$ for which the system is nonstiff.
In Figure~\ref{fig:fig1} we see numerical solutions computed by MPRK22($-0.5$) for two slightly different values of $\delta$. The rather large time step size $\Delta t = 1$, which is of course unsuitable for accuracy, was chosen to demonstrate the stability behavior of the scheme.
As described in Section~\ref{sec:alpha_neg} the MPRK22($-0.5$) scheme is stable and the iterates tend to the steady state of the corresponding initial value problem for all time step sizes, as long as the initial value is close enough to the steady state. This behavior is depicted in Figure~\ref{fig:fig1a}. 
Now, in Figure~\ref{fig:fig1b} the value of $\delta$, i.\,e.\ the distance to the steady state, is slightly increased and as a result the numerical approximations tend to a spurious steady state. As discussed above, we would like to have a linear stability behavior similar to A-stability, i.\,e.\ the stability should not depend on the initial values.
But here it is demonstrated that MPRK schemes exist for which stability crucially depends on the initial value.

\begin{figure}
\centering
\begin{subfigure}{0.5\textwidth}
    \includegraphics[width=\textwidth]{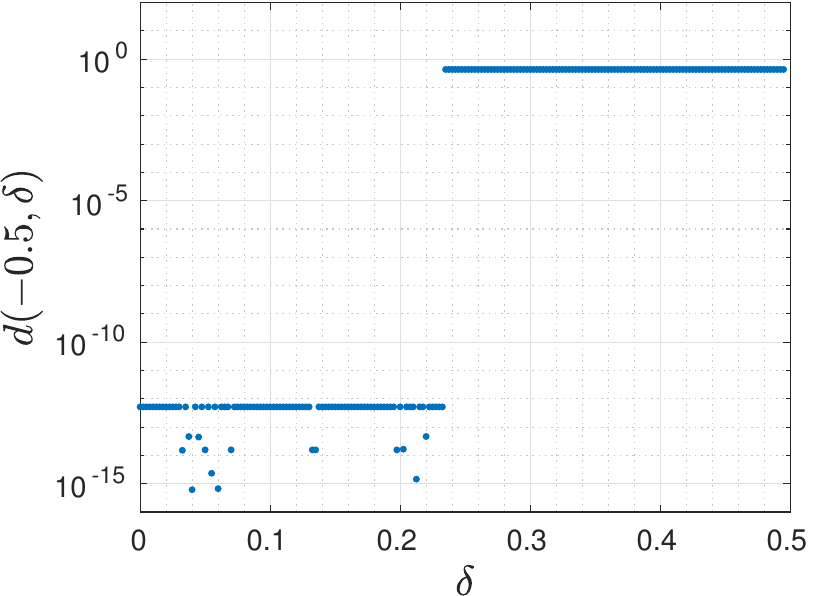}
    \caption{$a = 20$}
    \label{fig:fig2a}
\end{subfigure}%
\begin{subfigure}{0.5\textwidth}
    \includegraphics[width=\textwidth]{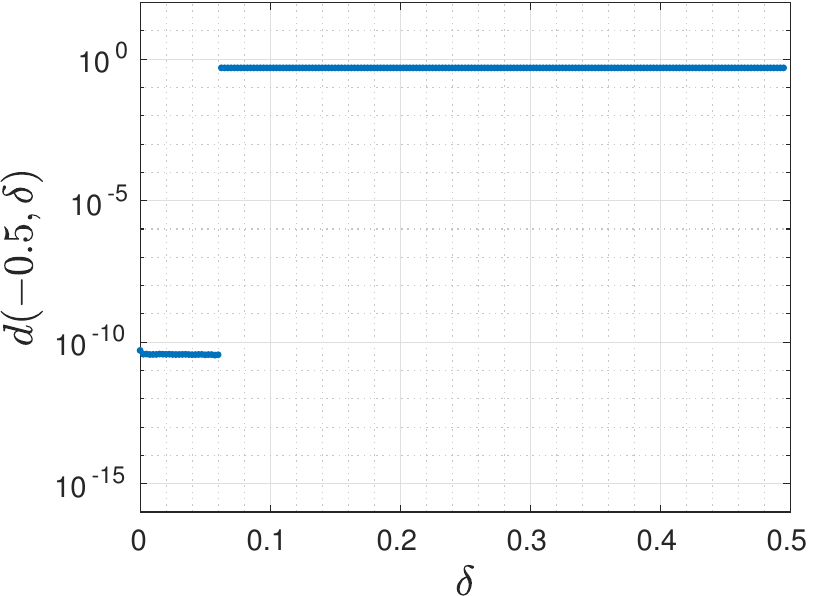}
    \caption{$a=200$}
    \label{fig:fig2b}
\end{subfigure}
\caption{Plots of $d(-0.5,\delta) = \norm{\vec y^* - \vec y^M}_{\infty}$ with 200 samples of $\delta$. The steady state of \eqref{eq:testprob} is $\vec y^*=(0.5,0.5)^T$ and $\vec y^M$ is computed by $M=10^4$ steps of MPRK22($-0.5$) with $\Delta t=1$.}
\label{fig:fig2}
\end{figure}

To gain a better insight into the dependence of the stability on $\delta$, we compute the steady states for $N=200$ equidistantly spaced samples of $\delta$ between $0$ and $0.5$. In each case the steady state is computed by performing $M=10^4$ steps of the MPRK22($-0.5$) scheme to obtain $\vec y^{M}$. To check if the correct or a spurious steady state is approached we compute the distance \[d(\alpha,\delta) = \norm{\vec y^* - \vec y^M}_{\infty}.\]
We expect $d(\alpha,\delta)\ll 1$ if the scheme is stable and tends to the steady state of the corresponding initial value problem and $d(\alpha,\delta)\gg 0$ if a spurious steady state is reached.
Figure~\ref{fig:fig2a} shows a plot of $d(-0.5,\delta)$ for $0\leq \delta < 0.5$. 
We see that the stability behavior changes abruptly from stable to unstable for some $\delta$ between $0.23$ and $0.24$ as suggested by Figure~\ref{fig:fig1}.
If we increase the stiffness of the system by choosing $a=200$ the critical value of $\delta$, for which the change from stable to unstable behavior occurs, decreases, as can be observed from Figure~\ref{fig:fig2b}. In this case the stability behavior changes from stable to unstable for $\delta\approx 0.06$.

\begin{figure}
\centering
\begin{subfigure}{0.5\textwidth}
    \includegraphics[width=\textwidth]{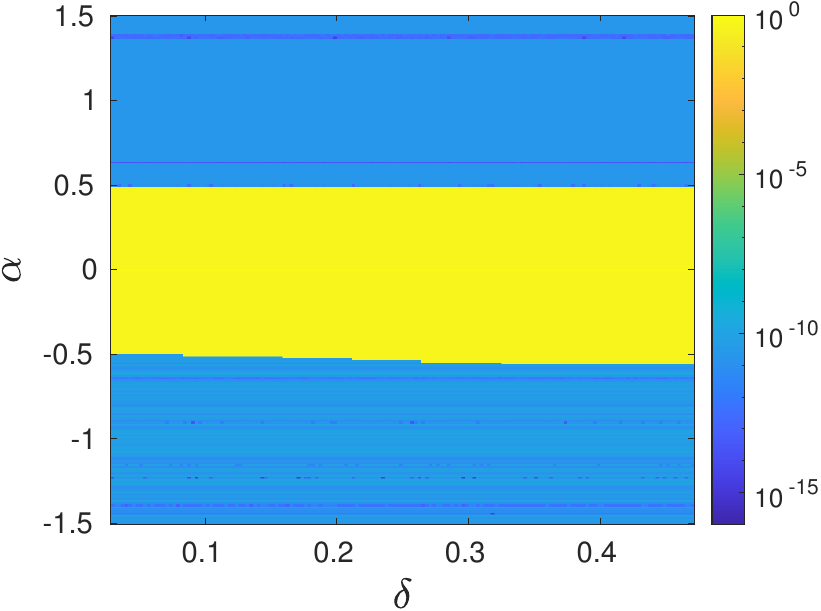}
    \caption{}
    \label{fig:fig3a}
\end{subfigure}%
\begin{subfigure}{0.5\textwidth}
    \includegraphics[width=\textwidth]{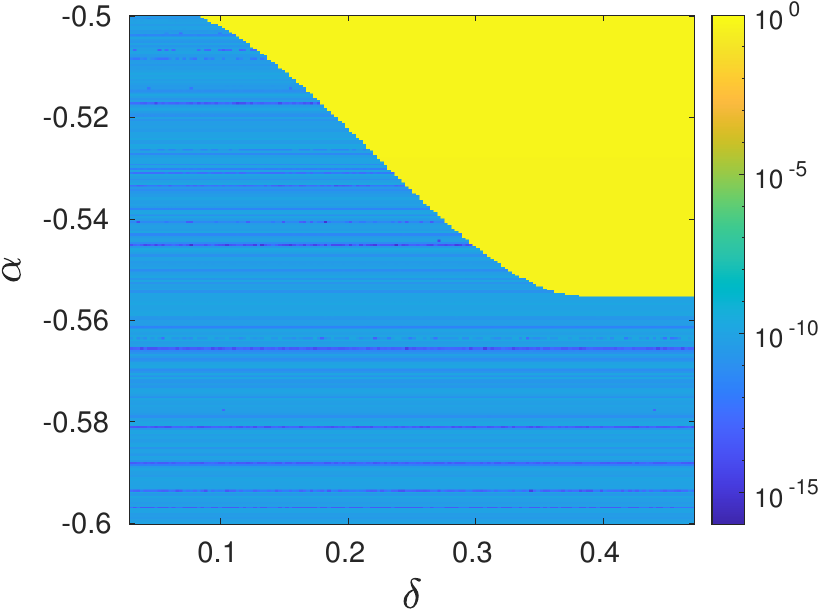}
    \caption{}
    \label{fig:fig3b}
\end{subfigure}
\caption{Plots of $d(\alpha,\delta) = \norm{\vec y^* - \vec y^M}_{\infty}$ with 160 samples of $\alpha$ and 241 samples $\delta$. The steady state of \eqref{eq:testprob} is $\vec y^*=0.5(1,1)^T$ and $\vec y^M$ is computed by $M=10^4$ steps of MPRK22($\alpha$) with $\Delta t=1$.}
\label{fig:fig3}
\end{figure}

Next, we want to answer the question if $\alpha=-0.5$ is the only value of $\alpha$ for which the linear stability of MPRK22($\alpha$) crucially depends on the initial value.  To answer this question, we again compute the distance $d(\alpha,\delta)$, but this time we vary both $\delta$ and $\alpha$. We use an equidistant grid with 160 samples of $\delta$ and 241 samples of $\alpha$ for $0< \delta < 0.5$ and $0<\abs{\alpha}\leq 2$. Again, we use $a=200$ in \eqref{eq:testprob}, $\Delta t=1$ and $M=10^4$ steps to compute the steady state. The result of this computation can be seen in Figure~\ref{fig:fig3a} and
we note that the MPRK22($\alpha$) scheme is not defined for $\alpha=0$, which is not indicated in the plot.
Based on Figure~\ref{fig:fig3a}, we can make the following statements.
First, the expected unstable behavior for $0<\abs{\alpha}<0.5$, which was discussed in Sections~\ref{sec:0<alpha<0.5} and \ref{sec:alpha_neg}, is clearly visible.  Second, for $\alpha\geq0.5$ we observe stable behavior independent of $\delta$. Therefore, local stability could actually be global stability for $\alpha\geq0.5$. Third, for $\alpha\leq -0.5$ there exists an  $\alpha^*<-0.5$ such that stability of MPRK22($\alpha$) is indeed only a local stability and for $\alpha<\alpha^*$ no unstable behavior can be observed. To see this more clearly, Figure~\ref{fig:fig3b} shows a zoom of Figure~\ref{fig:fig3a} with higher resolution for $-0.6\leq \alpha\leq -0.5$. As we can see, $\alpha=-0.5$ is not the only value of $\alpha$ for which the stability depends on the initial value. The same is true for all $\alpha$ between $\alpha^*\approx -0.56$ and $-0.5$.
 
\section{Summary and outlook}
We have discussed the linear stability behavior of MPRK22($\alpha$) schemes for all $\alpha\ne 0$ based on the local stability theory of \cite{IKM2022,IKM2022b}. 
Moreover, we have confirmed that negative RK parameters in MPRK schemes can lead to linear stability actually being a local property.
In particular, we have shown that there are MPRK22($\alpha$) schemes with $\alpha<-0.5$ that are linearly stable but converge to a spurious fixed point instead of the steady state of the corresponding linear system if the initial value is not sufficiently close to the steady state. 
Such a linear stability behavior is unknown from general linear methods. 
We further note that we cannot recommend to use MPRK22($\alpha$) schemes with $\alpha<0.5$ since apart from the time step size restrictions for $0<\alpha<0.5$, the stability of these numerical methods shows an artificial dependence on $\mat A^T$ and $\vec y^*$ that is not present in the continuous case. In addition, the convergence towards the steady state of the corresponding initial value problem cannot be guaranteed in the case of general linear systems by means of the stability theory from \cite{IKM2022,IKM2022b}, see \cite{OeffnerIzgin2023}.

Nevertheless, for MPRK22($\alpha$) schemes relevant in practice, i.\,e.\ $\alpha\approx 1$, the parameters of the underlying RK scheme are nonnegative and no restrictions on the initial value could be found numerically. 
Therefore, for these methods, it still seems that the local linear stability is actually a global linear stability. 

Since we would like to have MPRK schemes with properties similar to A-stable general linear methods, it is therefore essential that the convergence to the steady state of the corresponding initial value problem is given for all initial values which satisfy $\vec n_i^T\vec y^0=\vec n_i^T\vec y^*$, where the vectors $\vec n_i$ define the independent linear invariants of the system matrix. For this reason, there is a great need for a theory that makes statements not only about local stability but also about global stability. 

\begin{appendices}
\section{ }\label{apndx}
Here we show the computations of \cite{Schilling2023} that have been omitted in Section~\ref{sec:alpha_neg}. 

The MPRK22($\alpha$) scheme \eqref{eq:MPRK22<0} can be rewritten as
\begin{equation*}
\vec \Psi(\vec y^n, \vec y^{(2)}) = \vec 0,\quad
\vec \Phi(\vec y^n, \vec y^{(2)}, \vec y^{n+1}) = \vec 0,
\end{equation*}
with
\begin{subequations}\label{eq:PsiPhi}
\begin{align}
\vec \Psi(\vec u,\vec v)&=\vec u + a_{21}\Delta t\bigl(\diag(\vec v/\vec u)\mathbf B\vec u+\diag(\mathbf A)\diag(\vec v/\vec u)\vec u\bigr)-\vec v,\\
\vec \Phi(\vec u,\vec v,\vec w)&\!\begin{multlined}[t][10cm]
     = \vec u+\Delta t\bigl(b_1\mathbf B\diag(\vec w/\boldsymbol \sigma)\mathbf u+b_2\diag(\vec w/\boldsymbol \sigma)\mathbf B\vec v\\+b_1\diag(\mathbf w/\boldsymbol \sigma)\diag(\mathbf A)\mathbf u+b_2\diag(\mathbf A)\diag(\mathbf w/\boldsymbol \sigma)\mathbf v\bigr)-\vec w.
     \end{multlined}
\end{align}
\end{subequations}
This form is particularly helpful to compute the Jacobian $\mat D\vec g(\vec y^*)$ through implicit differentiation. To keep the notation short we denote by
$\mat D_{\vec x_j}\vec f(\vec x_1,\dots,\vec x_k)$ the Jacobian of a function $\vec f$ with respect to $\vec x_j$.
Furthermore, we define $\mat D^*_{\vec u}\Psi$ as $\mat D_{\vec u}\Psi(\vec y^*,\vec y^*)$ and $\mat D^*_{\vec u}\Phi$ as $\mat D_{\vec u}\Phi(\vec y^*,\vec y^*,\vec y^*)$. With this notation we have
\begin{equation}\label{eq:Dg}
\mat D\vec g(\vec y^*)=-(\mat D_{\vec w}^*\Phi)^{-1}\bigl(\mat D_{\vec u}^*\Phi-\mat D^*_{\vec v}\Phi(\mat D^*_{\vec v}\Psi)^{-1}\mat D^*_{\vec u}\Psi\bigr),\end{equation}
if $\Psi,\Phi\in\mathcal C^2$ and the inverse matrices $(\mat D_{\vec w}^*\Phi)^{-1}$ and $(\mat D^*_{\vec v}\Psi)^{-1}$ exist, see \cite{OeffnerIzgin2023}.

To actually compute the Jacobians appearing in \eqref{eq:Dg}, we rewrite \eqref{eq:PsiPhi} elementwise yielding
\begin{subequations}
\begin{align}
\Psi_i(\vec u,\vec v) &= -v_i+u_i - \alpha\Delta t\Biggl(\sum_{\substack{j=1 \\ j \neq i}}^N a_{ji}u_i v_j u_j^{-1} - \sum_{\substack{j=1 \\ j \neq i}}^N a_{ij}u_j v_i u_i^{-1}\Biggr),\label{eq:Psi_i}\\
\Phi_i(\vec u,\vec v,\vec w)&\!\begin{multlined}[t][10cm]
= -w_i + u_i + \Delta t \Biggl( 
\left(1-\tfrac{1}{2\alpha}\right) \sum_{\substack{j=1 \\ j \neq i}}^N a_{ij} w_j v_j^{-1/\alpha} u_j^{1/\alpha}\\
- \tfrac{1}{2\alpha} \sum_{\substack{j=1 \\ j \neq i}}^N a_{ji}v_i w_j v_j^{-1/\alpha} u_j^{-1+1/\alpha} 
 - \left(1-\tfrac{1}{2\alpha}\right) \sum_{\substack{j=1 \\ j \neq i}}^N a_{ji} w_i v_i^{-1/\alpha}u_i^{1/\alpha}\\
+ \tfrac{1}{2\alpha} \sum_{\substack{j=1 \\ j \neq i}}^N a_{ij} v_jw_iv_i^{-1/\alpha}u_i^{-1+1/\alpha} \Biggr).\label{eq:Phi_i}
\end{multlined}
\end{align}
\end{subequations}
By differentiation of \eqref{eq:Psi_i} we obtain the Jacobians 
\begin{align*}
(\mat D_{\vec u}\vec \Psi(\vec u,\vec v))_{iq} 
&=\begin{cases} 1-\alpha \Delta t \left(\sum_{\substack{j=1 \\ j \neq i}}^N  a_{ji} v_ju_j^{-1}+\sum_{\substack{j=1 \\ j \neq i}}^N  a_{ij} u_j v_iu_i^{-2}\right), & i=q,\\
-\alpha \Delta t \left(-a_{qi} u_i v_qu_q^{-2}-a_{iq} v_iu_i^{-1}\right), & i\ne q,
\end{cases}
\intertext{and}
(\mat D_{\vec v}\vec \Psi(\vec u,\vec v))_{iq} 
&=\begin{cases}  
-1+\alpha \Delta t \sum_{\substack{j=1 \\ j \neq i}}^N  a_{ij} u_j u_i^{-1},& i=q,\\
-\alpha \Delta t a_{qi} u_i u_q^{-1}, & i\ne q.
\end{cases}
\end{align*}
By differentiation of \eqref{eq:Phi_i} one can write
\[
(\mat D_{\vec u}\vec \Phi(\vec u,\vec v,\vec w))_{ii}  =  \begin{multlined}[t][12cm] + \Delta t \biggl(
-\tfrac{1}{\alpha} (1-\tfrac{1}{2\alpha}) 
\sum_{\substack{j=1 \\ j \neq i}}^N a_{ji} w_i v_i^{-1/\alpha} u_i^{-1+1/\alpha} \\[-\baselineskip]+ \left(-1+\tfrac{1}{\alpha}\right)  \tfrac{1}{2\alpha} 
\sum_{\substack{j=1 \\ j \neq i}}^N a_{ij} v_j w_iv_i^{-1/\alpha}u_i^{-2+1/\alpha}\biggr)\end{multlined}
\]
and
\[
(\mat D_{\vec u}\vec \Phi(\vec u,\vec v,\vec w))_{iq}=
\Delta t \Bigl(
\tfrac{1}{\alpha} \left(1-\tfrac{1}{2\alpha}\right) 
a_{iq} w_q v_q ^{-1/\alpha}u_q^{-1+1/\alpha}
-\left(-1+\tfrac{1}{\alpha}\right) \tfrac{1}{2\alpha} 
a_{qi} v_i w_q v_q ^{-1/\alpha} u_q^{-2+1/\alpha}\Bigr)
\]
for $i\ne q$.
Furthermore, we see
\[(\mat D_{\vec v}\vec \Phi(\vec u,\vec v,\vec w))_{ii}=\begin{multlined}[t][12cm]
\Delta t\biggr(
-\tfrac{1}{2\alpha} 
\sum_{\substack{j=1 \\ j \neq i}}^N a_{ji} w_j v_j^{-1/\alpha} u_j^{-1+1/\alpha}
+ \tfrac{1}{\alpha} \left(1-\tfrac{1}{2\alpha}\right) 
\sum_{\substack{j=1 \\ j \neq i}}^N  a_{ji} w_iv_i^{-1-1/\alpha}u_i^{1/\alpha}\\ 
 -\tfrac{1}{\alpha} \tfrac{1}{2\alpha}  
\sum_{\substack{j=1 \\ j \neq i}}^N  a_{ij} v_j w_iv_i^{-1-1/\alpha}u_i^{-1+1/\alpha}\biggl)
\end{multlined}\]
and
\[(\mat D_{\vec v}\vec \Phi(\vec u,\vec v,\vec w))_{iq}=\begin{multlined}[t][12cm] \Delta t \Bigl(
-\tfrac{1}{\alpha} \left(1-\tfrac{1}{2\alpha}\right) 
a_{iq} w_qv_q^{-1-1/\alpha}u_q^{1/\alpha} 
+ \tfrac{1}{\alpha} \tfrac{1}{2\alpha} a_{qi} v_i w_qv_q^{-1-1/\alpha}u_q^{-1+1/\alpha}\\ +\tfrac{1}{2\alpha}a_{iq}w_i v_i^{-1/\alpha} u_i^{-1+1/\alpha}\Bigr)\end{multlined}\]
for $i\ne q$.
Finally,
\[(\mat D_{\vec w}\vec \Phi(\vec u,\vec v,\vec w))_{ii}=-1+ \Delta t\biggl(
- \left(1-\tfrac{1}{2\alpha}\right)
\sum_{\substack{j=1 \\ j \neq i}}^N a_{ji} v_i^{-1/\alpha}u_i^{1/\alpha} 
+\tfrac{1}{2\alpha}
\sum_{\substack{j=1 \\ j \neq i}}^N a_{ij} v_j v_i^{-1/\alpha}u_i^{-1+1/\alpha}\biggr)\]
and
\[(\mat D_{\vec w}\vec \Phi(\vec u,\vec v,\vec w))_{iq}=\Delta t \Bigl(\left(1-\tfrac{1}{2\alpha}\right)
 a_{iq} v_q^{-1/\alpha}u_q^{1/\alpha}
-\tfrac{1}{2\alpha} a_{qi} v_i v_q^{-1/\alpha}u_q^{-1+1/\alpha}\Bigr)\]
hold for $i\ne q$.
Now, substituting $\vec u=\vec v=\vec w=\vec y^*$ yields the Jacobians
\begin{subequations}\label{eq:Jacobians}
\begin{align} 
\mat D_{\vec u}^*\vec\Psi
&= \mat I + \alpha \Delta t (\mat A + \diag(\vec y^{*}) \mat A^T \diag(\vec y^{*})^{-1}), \\
\mat D_{\vec v}^*\vec\Psi 
&= -\mat I - \alpha \Delta t \diag(\vec y^{*}) \mat A^T \diag(\vec y^{*})^{-1},\\
\mat D_{\vec u}^*\vec\Phi
&= \mat I + \Delta t \bigl(
\left(\tfrac{1}{\alpha}-\tfrac{1}{2\alpha^2}\right) \mat A
- \left(-\tfrac{1}{2\alpha} +\tfrac{1}{2\alpha^2}\right) \diag(\vec y^{*}) 
\mat A^T \diag(\vec y^{*})^{-1} \bigr) ,\\
\mat D_{\vec v}^*\vec\Phi 
&= \Delta t \bigl(
\left(-\tfrac{1}{2\alpha} +\tfrac{1}{2\alpha^2}\right) \mat A
+\tfrac{1}{2\alpha^2} \diag(\vec y^{*}) \mat A^T \diag(\vec y^{*})^{-1} \bigr),\\
\mat D_{\vec w}^*\vec\Phi 
&= -\mat I +\Delta t\bigl(
\left(1-\tfrac{1}{2\alpha}\right) \mat A
-\tfrac{1}{2\alpha} \diag(\vec y^{*}) \mat A^T \diag(\vec y^{*})^{-1} \bigr),
\end{align}
\end{subequations}
which shows that $\mat D\vec g(\vec y^*)$ as given in \eqref{eq:Dg} in general depends on $\mat A^T$ and $\vec y^*$. 
However, hereafter we only consider the linear system \eqref{eq:lin2x2} for which $\mathbf y^* = s(b,a)^T$ for some $s\in\R$ and consequently
\begin{multline}\label{eq:condition}\diag(\vec y^{*}) \mat A^T \diag(\vec y^{*})^{-1}=s\Mat{b & 0\\0 & a}\Mat{-a & a\\b & -b}\left(s^{-1}\Mat{b^{-1} & 0\\0 & a^{-1}}\right)\\=\Mat{-bab^{-1} & ba a^{-1}\\ a b b^{-1} & -aba^{-1}}=\mat A.\end{multline}
With this we can conclude
\begin{align*} 
\mat D_{\vec u}^*\vec\Psi
&= \mat I+2\alpha \Delta t \mat A, &
\mat D_{\vec v}^*\vec\Psi 
&= -\mat I - \alpha \Delta t  \mat A
\end{align*}
and
\begin{align*}
\mat D_{\vec u}^*\vec\Phi&=\mat I + \Delta t \bigl( \tfrac{3}{2\alpha} - \tfrac{1}{\alpha^2} \bigr) \mat A,&
\mat D_{\vec v}^*\vec\Phi &= \Delta t \bigl( -\tfrac{1}{2\alpha} + \tfrac{1}{\alpha^2} \bigr) \mat A, &
\mat D_{\vec w}^*\vec\Phi &= -\mat I + \Delta t \bigl( 1 - \tfrac{1}{\alpha} \bigr) \mat A.
\end{align*}
Inserting this into \eqref{eq:Dg} finally yields
\begin{multline}\label{eg:Dgy*}
\mat D\vec g(\vec y^*)
= -\bigl(-\mat I + \Delta t \bigl( 1 - \tfrac{1}{\alpha} \bigr) \mat A \big )^{-1}
\bigl(\mat I + \Delta t \bigl( \tfrac{3}{2\alpha} - \tfrac{1}{\alpha^2} \bigr) \mat A \\
-\Delta t \bigl( -\tfrac{1}{2\alpha} + \tfrac{1}{\alpha^2} \bigr) \mat A
\bigl(-\mat I-\alpha \Delta t \mat A \bigr)^{-1} \bigl( \mat I+2\alpha \Delta t \mat A \bigr) \bigr),
\end{multline}
under the assumption that $(\mat D_{\vec w}^*\Phi)^{-1}$ and $(\mat D^*_{\vec v}\Psi)^{-1}$ exist. 
To verify this we note that the system under consideration has eigenvalues $\lambda_1=0$ and $\lambda_2=-(a+b)<0$. Hence, if $\lambda\leq 0$ is an eigenvalue of $\mat A$, then $-1-\alpha\Delta t\lambda\leq -1$ is an eigenvalue of $\mat D^*_{\vec v}\Psi$, which shows that $\mat D^*_{\vec v}\Psi$ is invertible. Similarly, the eigenvalues of $\mat D_{\vec w}^*\Phi$ are $-1+\Delta t(1-\tfrac1{\alpha})\lambda\leq -1$. Thus, $\mat D_{\vec w}^*\Phi$ is invertible as well. 

The eigenvalues of \eqref{eg:Dgy*} are $R(\Delta t\lambda)$ with $\lambda$ being an eigenvalue of $\mat A$ and $R$ representing the stability function
\begin{equation*}
R(z)=
-\frac{1+\bigl( \tfrac{3}{2\alpha} - \frac{1}{\alpha^2} \bigr)z
-\bigl( -\tfrac{1}{2\alpha} + \frac{1}{\alpha^2} \bigr)z
\frac{1}{(-1-\alpha z)} \bigl(1+2\alpha z\bigr)}
{-1+\bigl(1-\frac{1}{\alpha}\bigr)z}.
\end{equation*}
We can rewrite $R$ as
\begin{equation*}
R(z)=
\frac{-(\alpha+2)z^2-(2\alpha^2+2)z-2\alpha}
{(2\alpha^2-2\alpha)z^2+(-2\alpha^2+2\alpha-2)z-2\alpha}
\end{equation*}
and compute its derivative
\begin{equation*}
R'(z)=\frac{2 \alpha^4 z^2 - \alpha^3 z^2 + 4 \alpha^3 z + 3 \alpha^2 z^2 - 2\alpha^2 z -3\alpha z^2
+2 \alpha^2 + 4\alpha z  + 2 z^2}
{2(\alpha z +1 )^2(\alpha z-\alpha-z)^2}.
\end{equation*}
For $z < 0$ and $\alpha < 0$ every term in the nominator is positive and the denominator is positive as well. Hence, $R'(z)>0$ for all $z < 0$ and $\alpha<0$, which implies that $R$ is monotonically increasing on $(-\infty, 0)$.
Furthermore, we find $R(0)=1$ and 
\begin{equation*}
\lim_{z \to -\infty} R(z)=- \frac{\alpha+2}{2\alpha(\alpha-1)}.
\end{equation*}
Together with the monotonicity of $R$ this results in
\begin{equation*}
-\frac{\alpha+2}{2\alpha(\alpha-1)}<R(z)<1 \quad \text{for } z<0 \text{ and } \alpha<0.
\end{equation*}
Since
\[-1\leq -\frac{\alpha+2}{2\alpha(\alpha-1)} \iff \alpha\leq-\frac12,\]
we can conclude
\begin{equation*}
|R(z)|<1 \quad \text{for } z<0 \text{ and } \alpha \leq -0.5.
\end{equation*}
According to the stability theory of \cite{IKM2022,IKM2022b} it follows that MPRK22($\alpha$) schemes with $\alpha \leq -0.5$ are unconditionally stable and converge to the steady state of the corresponding initial value problem  independent of the time step size $\Delta t$, if the initial value is close enough to the fixed point.
On the other hand, for $-0.5<\alpha<0$ there must exist $z<0$ with $\abs{R(z)}>1$. Hence, MPRK22($\alpha$) schemes with $-0.5<\alpha<0$ are only conditionally stable. Solving $R(z^*)=-1$ for $z^*$ shows
\[
z^*= \frac{2 \alpha^2 - \alpha +2+  \sqrt{4 \alpha^4+ 4 \alpha^3 - 3 \alpha^2 - 12 \alpha +4}}
{2 \alpha^2 - 3 \alpha -2}.
\]
Hence, stability for $-0.5<\alpha<0$ requires $z^*<z<0$.

\end{appendices}

\section*{Declarations}
\subsection*{Ethical Approval}
Not Applicable
\subsection*{Availability of supporting data}
The MATLAB code which was used to generate the numerical results is available from \texttt{https://github.com/SKopecz/locstabMPRK22.git}.
\subsection*{Competing interests}
The authors declare no competing interests.
\subsection*{Funding}
The author T.\ Izgin gratefully acknowledges the financial support by the Deutsche Forschungsgemeinschaft (DFG) through grant ME 1889/10-1.
\subsection*{Authors' contributions}
All authors have conceptualized the work and revised it critically. S.K. wrote the main manuscript text and performed the numerical experiments. A.S. performed the stability investigation in the appendix. T.I. wrote several paragraphs of the manuscript.
\subsection*{Acknowledgments}
Not Applicable 

\bibliography{locstabMPRK22}

\end{document}